\documentclass[12pt]{article}
\usepackage{amsfonts}
\input{epsf.sty}
\usepackage[dvips]{graphicx}
\usepackage{latexsym,amsmath,amsfonts,amscd, amsthm}
\usepackage{epsfig}
\usepackage{changebar}
\usepackage{pstricks}
\usepackage{pst-plot}
\usepackage{multirow,bigstrut}
\usepackage{subfigure}
\usepackage{placeins}
\usepackage{amssymb}
\usepackage{parskip}

\usepackage{mathrsfs}
\usepackage{amssymb,amsthm,hyperref,amsmath,graphicx,bm,color,booktabs}
\usepackage{bm}
\usepackage{comment}
\usepackage{enumerate}

\numberwithin{equation}{section} 

\newtheorem{thm}{Theorem}[section]

\newtheorem{lem}{Lemma}[section]
\newtheorem{cor}{Corollary}[section]

\theoremstyle{definition}

\newcommand{\brac}[1]{\left(#1\right)}
\newcommand{\abs}[1]{\left\vert#1\right\vert}
\newcommand{\norm}[1]{\left\Vert#1\right\Vert}
\def\half{\frac 1 2}

\newcommand{\diag}{\mbox{diag}}

\newcommand{\dt}{\Delta t}

\allowdisplaybreaks

\begin{document}

\baselineskip=2pc

\begin{center}
{\bf \large Uniform accuracy of implicit-explicit backward differentiation formulas (IMEX-BDF) for linear hyperbolic relaxation systems}
\end{center}

\vspace{.2in}

\centerline{
Zhiting Ma \footnote{Yanqi Lake Beijing Institute of Mathematical Sciences and Applications, Beijing 101408, China. E-mail: mazt@bimsa.cn}
\qquad
Juntao Huang \footnote{Department of Mathematics and Statistics, Texas Tech University, Lubbock, TX, 79409, USA. E-mail: juntao.huang@ttu.edu. Research is partially supported by NSF DMS-2309655 and DOE DE-SC0023164.}
\qquad
Wen-An Yong \footnote{Department of Mathematical Sciences, Tsinghua University, Beijing 100084, China; and Yanqi Lake Beijing Institute of Mathematical Sciences and Applications, Beijing 101408, China.
E-mail: wayong@tsinghua.edu.cn. Research is partially supported by National Key Research and Development Program of China (Grant no. 2021YFA0719200) and the National Natural Science Foundation of China (Grant no.12071246).
}
}

\vspace{.2in}

\centerline{\bf Abstract}
This work is concerned with the uniform accuracy of implicit-explicit backward differentiation formulas for general linear hyperbolic relaxation systems satisfying the structural stability condition proposed previously by the third author. 
We prove the uniform stability and accuracy of a class of IMEX-BDF schemes discretized spatially by a Fourier spectral method. The result reveals that the accuracy of the fully discretized schemes is independent of the relaxation time in all regimes. It is verified by numerical experiments on several applications to traffic flows, rarefied gas dynamics and kinetic theory.

\vfill

\noindent {\bf Keywords}: {}


\section{Introduction}
\setcounter{equation}{0}
\setcounter{figure}{0}
\setcounter{table}{0}

This paper is concerned with uniformly stable and accurate numerical methods for one-dimensional linear hyperbolic relaxation systems
\begin{equation}\label{eq:PDE-general}
	U_t + A U_{x} = \frac{1}{\varepsilon}Q U.
\end{equation}
Here $U=U(x,t)\in\mathbb{R}^n$, $x \in\mathbb{R}$, $t\geq 0$, $A $ and $Q$ are $n\times n$ constant matrices, the subscripts $t$ and $x$ refer to the partial derivatives with respect to $t$ and $x$, and $\varepsilon > 0$ is a small parameter standing for the relaxation time. 
Such partial differential equations (PDEs) are the linearized version of first-order PDEs with relaxation. The latter models a large number of different irreversible phenomena. Important examples include kinetic theories (moment closure systems \cite{levermore1996moment,Di2017nm}, discrete-velocity kinetic models \cite{broadwell1964shock,platkowski1988discrete}), nonlinear optics \cite{hanouzet2000approximation}, radiation hydrodynamics \cite{pomraning2005equations,mihalas2013foundations}, traffic flows \cite{aw2000siam}, dissipative relativistic fluid flows \cite{geroch1990dissipative}, chemically reactive flows \cite{giovangigli2012multicomponent}, and invisicid gas dynamics with relaxation \cite{zeng1999gas}.

Due to the small parameter $\varepsilon$, usual numerical schemes are stable only if the time step is of order $O(\varepsilon)$. To overcome this drawback, the so-called implicit-explicit (IMEX) schemes were adapted \cite{ascher1995siam}, where the convection part is treated explicitly and the source term is treated implicitly. 
The IMEX schemes include the IMEX Runge-Kutta  method (IMEX-RK, e.g., \cite{ascher1997,kennedy2003additive,dimarco2013siam,pareschi2005jsc}) and IMEX multistep method including IMEX backward differentiation formulas (IMEX-BDF, e.g.,  \cite{ascher1995siam,hundsdorfe2007imex,dimarco2017siam,albi2020implicit-explicit}). 
As reported in \cite{boscarino2007error-analysis,hu_uniform_2019}, many IMEX-RK schemes suffer from accuracy degeneration when $\varepsilon$ goes to zero, while the numerical experiments indicate the uniform accuracy of certain IMEX multistep schemes for a wide range of $\varepsilon$ \cite{hundsdorfe2007imex,dimarco2017siam,albi2020implicit-explicit}.

The aim of this work is to clarify the uniform accuracy of the multistep schemes for linear hyperbolic relaxation systems \eqref{eq:PDE-general}.
For the Jin-Xin model \cite{jin1995cpam} as a specific relaxation system, the uniform stability and accuracy have been studied in \cite{hu_uniform_2019} for the IMEX-BDF schemes and in \cite{hu_uniform_2023} for the IMEX-RK schemes.
Thus, our task is to generalize the analysis in \cite{hu_uniform_2019} for the Jin-Xin model to general hyperbolic relaxation systems satisfying the structural stability condition proposed in \cite{yong_singular_1999}. As shown in \cite{yong_singular_1999,liu2001basic,Yong2008}, the structural stability condition are tacitly respected by many well-developed physical theories. Therefore, our analysis is expected to have a wide range of applications.  

Under the structural stability condition, we prove the uniform stability and accuracy of the fully discretized IMEX-BDF schemes up to fourth order. The spatial discretization is done by adopting a Fourier spectral method \cite{hesthaven2007spectral}. The proof invokes a multiplier technique developed in \cite{dahlquist1978G-stability,nevanlinna1981multiplier}. 
Our results hold for any value of the small parameter $\varepsilon$. In other words, the accuracy of the schemes is independent of $\varepsilon$ in all regimes.
We also present numerical tests to verify our theoretical results with several specific relaxation systems, including the linearized Aw-Rascle-Zhang traffic model \cite{aw2000siam,zhang2002non}, the Broadwell model \cite{broadwell1964shock}, and a moment closure system \cite{grad1949kinetic,cai2014cpam}.

The rest of the paper is organized as follows. In Section \ref{sec:preliminaries}, we introduce the structural stability condition and IMEX-BDF schemes for the relaxation systems \eqref{eq:PDE-general}. Section \ref{sec:2-constant} is devoted to our main results including uniform-in-$\varepsilon$ stability and accuracy of the IMEX-BDF schemes.  Numerical experiments are presented in Section \ref{sec:numerical-tests} to validate our theoretical findings.

\section{Preliminaries}\label{sec:preliminaries}
In this section, we introduce the structural stability condition and a class of implicit-explict backward differentiation formulas (IMEX-BDF) for system \eqref{eq:PDE-general}.

\subsection{Structural Stability Condition}\label{subsec:ssc}

The structural stability condition reads as
\begin{enumerate}[(i)]
	\item There is an invertible $n\times n$ matrix $P$ and an invertible $r\times r$ $(0<r\le n)$ matrix $\hat{S}$ such that
	\begin{equation}\nonumber
		P Q = 
		\left(
		\begin{array}{cc}
			0 & 0 \\
			0 & \hat{S}
		\end{array}
		\right)
		P.
	\end{equation}

    \item There exists a symmetric positive-definite (SPD) matrix $A_0$ such that
	\begin{equation}\nonumber
		A_0  A  = A^T  A_0.
	\end{equation}

	\item The hyperbolic part and the source term are coupled in the sense:
	\begin{equation}\nonumber
		A_0  Q  + Q^T  A_0  \le -P^T 
		\left(
		\begin{array}{cc}
			0 & 0 \\
			0 & I_r
		\end{array}
		\right)
		P .
	\end{equation}
\end{enumerate}
Here the superscript $T$ denotes the transpose and $I_r$ is the unit matrix of order $r$.

About this set of conditions, we remark as follows. 
Condition (i) is classical for initial-value problems of systems of ordinary differential equations (ODE, spatially homogeneous systems), while (ii) means the symmetrizable hyperbolicity of the system of first-order partial differential equations (PDE) in \eqref{eq:PDE-general}. Condition (iii) characterizes a kind of coupling between the ODE and PDE parts. 
As shown in \cite{yong_singular_1999,liu2001basic,Yong2008}, the structural stability condition has been tacitly respected by many well-developed physical theories. 
Recently, it is shown in \cite{Di2017nm,zhao2017stability,ma2023nonrelativisti} to be proper for certain moment closure systems. 
Under the structural stability condition, the existence and stability of the zero relaxation limit of the corresponding initial-value problems have been established in \cite{yong_singular_1999}. 

Assuming the structural stability condition, we introduce $\tilde{U} := PU$
and transform system \eqref{eq:PDE-general} into its equivalent version
\begin{equation}\nonumber
	\tilde{U}_t + \tilde{A} \tilde{U}_{x} = \frac{1}{\varepsilon}
	\left(
	\begin{array}{cc}
		0 & 0 \\
		0 & \hat{S}
	\end{array}
	\right)
	\tilde{U},
\end{equation}
where $\tilde{A} := P A P^{-1}$.
It is easy to see that the above equivalent version satisfies the structural stability condition with $\tilde{P}=I$ and $\tilde{A}_0 = P^{-T}A_0P^{-1}$. 
Thus, throughout this paper we only consider the transformed version (drop the tilde)
\begin{equation}\label{eq:PDE}
	U_t + AU_{x} = \frac{1}{\varepsilon}\left(
\begin{array}{cc}
	0 & 0 \\
	0 & \hat{S}
\end{array}
\right) U \equiv \frac{1}{\varepsilon} Q U .
\end{equation}

It was proved in \cite{yong_singular_1999} (Theorem 2.2) that $P^{-T}A_0P^{-1}$ is a block-diagonal matrix (with the same partition as in (i) and (iii)).
Thus, the symmetrizer for \eqref{eq:PDE} has the following block-diagonal form
\begin{equation}\nonumber
A_0 =
\left(
\begin{array}{cc}
	A_{01} & 0 \\
	0 & A_{02}
\end{array}
\right).
\end{equation}
We further assume that $A_{02}\hat{S}$ is symmetric (negative-definite), which holds true for many physical models \cite{Yong2008}.

\subsection{IMEX-BDF schemes}\label{subsec:IMEX-BDF}

Let $u^{n}=u^n(x)$ denote the numerical solution at time $t_n = T_0 + n\Delta t$, where $T_0$ is the initial time, $n$ is a non-negative integer, and $\Delta t$ is the time step. The $q$-th order IMEX-BDF scheme for system \eqref{eq:PDE}
read as
\begin{equation}\label{eq:BDF-ODE}
	\sum_{i=0}^q \alpha_i u^{n+i} + \dt \sum_{i=0}^{q-1} \gamma_i A  u_x^{n+i} = \beta \frac{\dt }{\varepsilon}Q u^{n+q}.
\end{equation}
Here, $\alpha:=(\alpha_0,\dots,\alpha_q)$, $\gamma:=(\gamma_0,\dots,\gamma_{q-1})$ and $\beta > 0$ are constant to be determined by the requirement that \eqref{eq:BDF-ODE} is $q$-th order accurate \cite{hundsdorfe2007imex}. Examples are 
\begin{itemize}
\item 	
$q=1$:
\begin{equation}
	u^{n+1} - u^n + \dt A u_x^n =  \frac{\dt}{\varepsilon} Q u^{n+1},
\end{equation}

\item 
$q=2$:
\begin{equation}\nonumber
	u^{n+2} - \frac{4}{3}u^{n+1} + \frac{1}{3}u^n + \dt \brac{ \frac{4}{3}A u_x^{n+1} - \frac{2}{3} A u_x^n } =  \frac{2}{3} \frac{\dt}{\varepsilon}Q u^{n+2},
\end{equation}

\item 
$q=3$:
\begin{equation}\nonumber
	u^{n+3} - \frac{18}{11} u^{n+2} + \frac{9}{11}u^{n+1} - \frac{2}{11}u^n + \dt \brac{\frac{18}{11}A u_x^{n+2} - \frac{18}{11}A u_x^{n+1} + \frac{6}{11 } A u_x^n} = \frac{6}{11}\frac{\dt}{\varepsilon}Q u^{n+3}.
\end{equation}

\end{itemize}

To analyze the IMEX-BDF scheme \eqref{eq:BDF-ODE}, we need the following multiplier technique established in  \cite{akrivis2016backward,hu_uniform_2019}.
\begin{lem}\label{lem:multiplier-Hu}
	Given $\alpha= (\alpha_0,\dots,\alpha_q)$ and $\gamma=(\gamma_0,\dots,\gamma_{q-1})$ in \eqref{eq:BDF-ODE} with $q=1, 2, 3, 4$, there exist a positive-definite quadratic form 
$$
G(u_1, \dots, u_q) = \sum_{i, j=1}^q g_{ij}u_i u_j,
$$
a semi-positive-definite quadratic form
$$
A(u_1, \dots, u_{q-1}) = \sum_{i, j=1}^{q-1} a_{ij} u_i u_j,
$$
two linear forms $L_1(u_1, \dots, u_{q-1})$,  $L_2(u_1, \dots, u_{q})$ such that
	\begin{equation}\nonumber
		\begin{aligned}
		    &\left(u_q - L_1(u_1, \dots, u_{q-1}\right))\sum_{i=0}^q\alpha_i u_i \\
      =& G(u_1, \dots, u_q) - G(u_0, \dots, u_{q-1}) + d_1\left(u_q - L_1(u_1, \dots, u_{q-1}) - d_2\sum_{i=0}^{q-1}\gamma_i u_i\right)^2
		\end{aligned}
	\end{equation}
	and
	\begin{equation}\nonumber
		\left(u_q - L_1(u_1, \dots, u_{q-1})\right)u_q = A(u_2, \dots, u_{q}) - A(u_1, \dots, u_{q-1}) + L_2^2(u_1, \dots, u_{q}).
	\end{equation}
  Here constants $d_1>0$ and all other constants are real.
\end{lem}

The proof of this lemma can be found in \cite{hu_uniform_2019}. Here we
 list the quadratic forms, the linear forms, and the constants for $q=1, 2$.
 \begin{itemize}
     \item $q=1$:
     \begin{equation}\nonumber
		G(u_1) = \frac{1}{2}u_1^2, \quad d_1 = \frac{1}{2}, \quad d_2 = 1, \quad L_2(u_1) = u_1.
	\end{equation}
 \item $q=2$:
  \begin{equation}\nonumber
		\begin{aligned}
		    &G(u_1, u_2) = \frac{1}{6}u_1^2 - \frac{2}{3}u_1u_2 + \frac{5}{6}u_2^2, \quad A(u_1) = 0, \quad L_1(u_1) = 0,\\[4mm]
             &L_2(u_1, u_2) = u_2, \quad d_1 = \frac{1}{6}, \quad d_2 = \frac{3}{2}.
		\end{aligned}
\end{equation}
 \end{itemize}

For our purpose, we generalize Lemma \ref{lem:multiplier-Hu} to the case where $u_j$  are vectors. 
To do this, we take a symmetric positive-definite (SPD) matrix $H$ and define a 
weighted inner-product for vectors $u,v\in\mathbb{R}^n$: 
\begin{equation}\nonumber
	(u, v)_{H} := u^T H v
\end{equation}
and norm
\begin{equation}\nonumber
	\norm{u}_{H} := \sqrt{(u,u)_H}.
\end{equation}
When $H=I_n$, the subscript $H$ will be omitted.

The generalized version of Lemma \ref{lem:multiplier-Hu} is 
\begin{lem}\label{lemma:multiplier}
Let the coefficients $g_{ij}$, $a_{ij}$ of the quadratic forms $G(u_1, \dots, u_q)$ and $A(u_1, \dots, u_{q-1})$, $L_1(u_1, \dots, u_{q-1}),  L_2(u_1, \dots, u_q)$ and $d_1, d_2$   be same as those in Lemma \ref{lem:multiplier-Hu}. 
For $u_j \in \mathbb{R}^n (j=0, \cdots, q)$ with $q=1,2,3,4$, set 
\begin{equation*}
   \begin{aligned}
        G_H(u_1, \dots, u_q) = \sum_{i,j=1}^q g_{ij}(u_i, u_j)_H,\quad A_H(u_1, \dots, u_{q-1}) = \sum_{i,j=1}^{q-1} a_{ij}(u_i, u_j)_H.
   \end{aligned}
\end{equation*}
Then the following two equalities hold:
	\begin{equation}\nonumber
		\begin{aligned}
		    &\left(u_q - L_1(u_1, \dots, u_{q-1}), \sum_{i=0}^q\alpha_i u_i\right)_H \\
      =& G_H(u_1, \dots, u_q) - G_H(u_0, \dots, u_{q-1}) + d_1 \norm{u_q - L_1(u_1, \dots, u_{q-1})  - d_2\sum_{i=0}^{q-1}\gamma_i u_i}_H^2\\
		\end{aligned}
	\end{equation}
	and
	\begin{equation}\nonumber
		\left(u_q - L_1(u_1, \dots, u_{q-1}), u_q\right)_H = A_H(u_2, \dots, u_{q}) - A_H(u_1, \dots, u_{q-1}) + \norm{L_2(u_1, \dots, u_q)}_H^2.
	\end{equation}
\end{lem}
\begin{proof}
It is well-known that for the given SPD matrix $H$, there exists a SPD matrix $M$ such that $H=M^2$. Then, for $u, v\in\mathbb{R}^n$ define $\tilde{u}=Mu$ and $\tilde{v}=Mv$. It holds that
	\begin{equation}\nonumber
		(u, v)_H = (\tilde{u}, \tilde{v}), \quad \norm{u}_H = \norm{\tilde{u}}.
	\end{equation}
 Thus the right-hand side of the first equality is equal to
	\begin{equation}\nonumber
		\begin{aligned}
			\textrm{RHS} 
			={}& \sum_{i,j=1}^q g_{ij}(u_i, u_j)_H - \sum_{i,j=0}^{q-1} g_{ij}(u_i, u_j)_H + d_1 \norm{u_q - L_1(u_1, \dots, u_{q-1}) - d_2 \sum_{i=0}^{q-1}\gamma_i u_i}_H^2\\
			={}& \sum_{i,j=1}^q g_{ij}(\tilde{u}_i, \tilde{u}_j) - \sum_{i,j=0}^{q-1} g_{ij}(\tilde{u}_i, \tilde{u}_j) + d_1 \norm{\tilde{u}_q - L_1(\tilde{u}_1, \dots, \tilde{u}_{q-1}) - d_2 \sum_{i=0}^{q-1}\gamma_i \tilde{u}_i}^2 \\
			={}& \left(\tilde{u}_q - L_1(\tilde{u}_1, \dots, \tilde{u}_{q-1}), ~\sum_{i=0}^q\alpha_i \tilde{u}_i\right) \\
			={}& \left(u_q - L_1(u_1, \dots, u_{q-1}), ~\sum_{i=0}^q\alpha_i u_i\right)_H = \textrm{LHS}.
		\end{aligned}
	\end{equation}
 Here the third equality follows from Lemma \ref{lem:multiplier-Hu} for each component of the $n$-vectors. Similarly, the second equality can be shown. This completes the proof.
\end{proof}

\section{Uniform accuracy} \label{sec:2-constant}

In this section, we consider system
\eqref{eq:PDE} 
with periodic boundary conditions. 
As in \cite{hu_uniform_2019}, 
we use the Fourier-Galerkin spectral method to the semi-discretized IMEX-BDF scheme \eqref{eq:BDF-ODE} in the spatial direction to obtain 
\begin{equation}\label{scheme1}
	\sum_{i=0}^q\alpha_i U_N^{n+i} + \dt A \sum_{i=0}^{q-1}\gamma_i (U_N^{n+i})_x = \frac{\beta\dt}{\varepsilon} QU_N^{n+q}.
\end{equation}
Here $U_N \in P_N := \textrm{span}\{ e^{ikx}| -N\le k\le N \}$ with $N$ being an integer.
For $P_N$-functions $U_N$, the following inequality is known \cite{hesthaven2007spectral}:  
\begin{equation}\label{equ:spectral}
	\norm{(U_N)_x}^2  \le N^2  \norm{U_N}^2.
\end{equation}

Here the notation $\norm{\cdot}$ denotes the usual $L^2$ norm of the square integrable periodic functions.

\subsection{Stability} \label{subsec:stability}
Assume the structural stability condition and the symmetry of the matrix $A_{02}\hat{S}$. In this subsection, we analyze the uniform-in-$\varepsilon$ stability of the fully discretized scheme \eqref{scheme1}. 
The main idea of our analysis will be illustrated firstly with the first-order scheme.

\subsubsection{First-order scheme}

For $q=1$, scheme \eqref{scheme1} reads as
\begin{equation}\nonumber
	U_N^{n+1} - U_N^{n} + \dt A (U_N^n)_x = \frac{\dt}{\varepsilon} QU_N^{n+1}.
\end{equation}
Multiplying this scheme with $(U_N^{n+1})^T A_0$ and integrating the resultant equality over $x$ gives 
\begin{equation}\label{eq:stablity-imex1}
	\int(U_N^{n+1})^T A_0(U_N^{n+1} - U_N^{n}) + \dt \int (U_N^{n+1})^T A_0 A (U_N^n)_x  = \frac{\dt}{\varepsilon} \int (U_N^{n+1})^T A_0 QU_N^{n+1} .
\end{equation}
Since $A_0$ is symmetric, the first term on the LHS of \eqref{eq:stablity-imex1} can be decomposed as
\begin{equation}\nonumber
	\begin{aligned}
		& \int(U_N^{n+1})^T A_0(U_N^{n+1} - U_N^{n}) \\
		={}& \frac{1}{2}\int(U_N^{n+1})^T A_0 U_N^{n+1}  - \frac{1}{2}\int(U_N^{n})^T A_0 U_N^{n}  + \frac{1}{2}\int(U_N^{n+1}-U_N^{n})^T A_0 (U_N^{n+1}-U_N^{n}),
	\end{aligned}
\end{equation}
while the second term is
\begin{equation}\nonumber
	\begin{aligned}
		& \dt \int (U_N^{n+1})^T A_0 A (U_N^n)_x  \\
		={}& \dt \int (U_N^{n+1}-U_N^n)^T A_0 A (U_N^n)_x  + \dt \int (U_N^{n})^T A_0 A (U_N^n)_x  \\
		={}& \dt \int (U_N^{n+1}-U_N^n)^T A_0 A (U_N^n)_x  + \frac{1}{2}\dt \int ((U_N^{n})^T A_0 A U_N^n)_x  \\
		={}& \dt \int (U_N^{n+1}-U_N^n)^T A_0 A (U_N^n)_x .
	\end{aligned}
\end{equation}
Here we have used the symmetry of $A_0A$ and the periodic boundary conditions. Thanks to the structural stability condition (iii), the RHS in \eqref{eq:stablity-imex1} is non-negative. Thus, it follows from \eqref{eq:stablity-imex1} that 
\begin{equation}\nonumber
	\begin{aligned}
		 \frac{1}{2}\int(U_N^{n+1})^T A_0 U_N^{n+1}  - \frac{1}{2}\int(U_N^{n})^T A_0 U_N^{n}  & + \frac{1}{2}\int(U_N^{n+1}-U_N^{n})^T A_0 (U_N^{n+1}-U_N^{n})  \\
		& + \dt \int (U_N^{n+1}-U_N^n)^T A_0 A (U_N^n)_x  \le 0.
	\end{aligned}
\end{equation}

Define $E^n := \frac{1}{2}\int(U_N^{n})^T A_0 U_N^{n} $ and denote by $2 \kappa$ the smallest eigenvalue of the SPD matrix $A_0$. We deduce from the last inequality and the inequality \eqref{equ:spectral} that
\begin{equation}\nonumber
	\begin{aligned}
		& E^{n+1} - E^n \\
		\le{}& -\frac{1}{2}\int(U_N^{n+1}-U_N^{n})^T A_0 (U_N^{n+1}-U_N^{n}) 
		 - \dt \int (U_N^{n+1}-U_N^n)^T A_0 A (U_N^n)_x  \\[4mm]
		\le{}& -\kappa \norm{U_N^{n+1}-U_N^n}^2 
		 + \kappa  \norm{U_N^{n+1}-U_N^n}^2  + \frac{C(\dt)^2}{\kappa} \norm{(U_N^n)_x}^2  \\[4mm]
  \le{}& \frac{C(\dt)^2}{\kappa}N^2  \norm{U_N^n}^2 \leq \frac{C(\dt)^2}{\kappa^2}N^2  E^n.
	\end{aligned}
\end{equation}

Finally, let $\dt\le c_{CFL}/N^2$. Then we have
\begin{equation}\nonumber
	E^{n} \le (1+C\dt) E^{n-1} \le (1+C\dt)^n E_0 \le e^{CT}E_0,
\end{equation}
namely,
\begin{equation}\nonumber
	\int(U_N^{n})^T A_0 U_N^{n}  \le e^{CT}\int(U_N^{0})^T A_0 U_N^{0} .
\end{equation}
This is the stability of the first-order fully discretized IMEX-BDF scheme.

\subsubsection{Higher-order schemes}

For other $q$, we have the following similar stability result.  
\begin{thm}\label{thm:stability-const}
	Under the structural stability condition, assume the CFL condition $\dt\le c_{CFL}/N^2$ with $c_{CFL}>0$ a constant. Then the IMEX-BDF scheme  \eqref{scheme1} with $q=1,2,3,4$ is uniformly stable in the sense that
	\begin{equation}\nonumber
		\norm{U_N^n}^2\le C\sum_{i=0}^{q-1}\brac{\norm{U_N^i}^2 + \frac{\dt}{\varepsilon}\norm{W_N^i}^2}
	\end{equation}
	for integer $n$ such that $t_n=T_0+n\dt\le T$, where $C$ is a constant independent of $\varepsilon$, $N$ and $\dt$, and
 $U^n_N=\begin{pmatrix}
      V^n_N\\W^n_N
  \end{pmatrix}$.
\end{thm}

\begin{proof}
Recall the scheme \eqref{scheme1} 
\begin{equation}\nonumber
	\sum_{i=0}^q\alpha_i U_N^{n+i} + \dt A \sum_{i=0}^{q-1}\gamma_i (U_N^{n+i})_x = \frac{\beta\dt}{\varepsilon} QU_N^{n+q}.
\end{equation} 
In Lemma \ref{lemma:multiplier}, taking $H=A_0$ from the structural stability condition we have
\begin{equation}\nonumber
	\begin{aligned}
		& \int \left(U_N^{n+q}-\sum_{i=1}^{q-1}\eta_i U_N^{n+i} \right)^T A_0 \sum_{i=0}^q\alpha_i U_N^{n+i} \\
		={}& \int G_{A_0}(U_N^{n+1}, \dots, U_N^{n+q}) - \int G_{A_0}(U_N^{n}, \dots, U_N^{n+q-1})\\
  {}&+ d_1 \norm{U_N^{n+q} - \sum_{i=1}^{q-1}\eta_i U_N^{n+i} - d_2 \sum_{i=0}^{q-1}\gamma_i U_N^{n+i}}_{A_0}^2.
	\end{aligned}
\end{equation}
Thanks to the symmetry of $A_0A$ and the periodic boundary conditions, we deduce from the inequality \eqref{equ:spectral} that
\begin{equation}\nonumber
	\begin{aligned}
		&\abs{ \dt \int \left(U_N^{n+q}-\sum_{i=1}^{q-1}\eta_i U_N^{n+i} \right)^T A_0 A \sum_{i=0}^{q-1}\gamma_i (U_N^{n+i})_x }\\
		\leq {}&\abs{ \dt \int \left(U_N^{n+q}-\sum_{i=1}^{q-1}\eta_i U_N^{n+i} - d_2 \sum_{i=0}^{q-1}\gamma_i U_N^{n+i} \right)^T A_0 A \sum_{i=0}^{q-1}\gamma_i (U_N^{n+i})_x }\\
		&+\abs{ \dt \int d_2\left( \sum_{i=0}^{q-1}\gamma_i U_N^{n+i} \right)^T A_0 A \sum_{i=0}^{q-1}\gamma_i (U_N^{n+i})_x }\\
		\leq {}& \abs{  \dt \int \left(U_N^{n+q}-\sum_{i=1}^{q-1}\eta_i U_N^{n+i} - d_2 \sum_{i=0}^{q-1}\gamma_i U_N^{n+i} \right)^T A_0 A \sum_{i=0}^{q-1}\gamma_i (U_N^{n+i})_x }\\
		\le{}& \kappa \norm{U_N^{n+q}-\sum_{i=1}^{q-1}\eta_i U_N^{n+i} - d_2 \sum_{i=0}^{q-1}\gamma_i U_N^{n+i}}^2  +  \frac{C(\dt)^2}{\kappa}  \norm{\sum_{i=0}^{q-1} \gamma_i (U_N^{n+i})_x}^2 \\
		\le{}& \kappa \norm{U_N^{n+q}-\sum_{i=1}^{q-1}\eta_i U_N^{n+i} - d_2 \sum_{i=0}^{q-1}\gamma_i U_N^{n+i}}^2  +  \frac{C(\dt)^2N^2}{\kappa}  \sum_{i=0}^{q-1} \norm{U_N^{n+i}}^2  \\
	\end{aligned}
\end{equation}
with $\kappa>0$. 
Moreover, 
the source term can be estimated as
\begin{equation}\nonumber
	\begin{aligned}
		& \int \left(U_N^{n+q}-\sum_{i=1}^{q-1}\eta_i U_N^{n+i} \right)^T A_0 \frac{\beta\dt}{\varepsilon} QU_N^{n+q}  \\
		={}& \frac{\beta\dt}{\varepsilon} \int \left(U_N^{n+q}-\sum_{i=1}^{q-1}\eta_i U_N^{n+i} \right)^T A_0  QU_N^{n+q}  \\
		={}& \frac{\beta\dt}{\varepsilon} \int \left(U_N^{n+q}-\sum_{i=1}^{q-1}\eta_i U_N^{n+i} \right)^T 
		\left(
		\begin{array}{cc}
			0 & 0 \\
			0 & A_{02}\hat{S}
		\end{array}
		\right) 
		U_N^{n+q}  \\
		={}& -\frac{\beta\dt}{\varepsilon} \int \left(W_N^{n+q}-\sum_{i=1}^{q-1}\eta_i W_N^{n+i} \right)^T M W_N^{n+q}  \\
		={}& - \frac{\beta\dt}{\varepsilon} \brac{ \int A_{M}(W_N^{n+2}, \dots, W_N^{n+q}) - \int A_{M}(W_N^{n+1}, \dots, W_N^{n+q-1}) + \norm{\sum_{i=1}^q c_i W_N^{n+i}}^2_M}
	\end{aligned}
\end{equation}
with $U=\begin{pmatrix}
    V\\ W
\end{pmatrix}$ and $M:= - A_{02}\hat{S}$ a SPD matrix.
Combining the last three estimates, we arrive at
\begin{equation}\label{eq:combine-estimate}
	\begin{aligned}
		& \int G_{A_0}(U_N^{n+1}, \dots, U_N^{n+q}) - \int G_{A_0}(U_N^{n}, \dots, U_N^{n+q-1})\\
    {}&+ d_1 \norm{U_N^{n+q} - \sum_{i=1}^{q-1}\eta_i U_N^{n+i} - d_2 \sum_{i=0}^{q-1}\gamma_i U_N^{n+i}}_{A_0}^2 \\[4mm]
		\le{} & \kappa  \norm{U_N^{n+q}-\sum_{i=1}^{q-1}\eta_i U_N^{n+i} - d_2 \sum_{i=0}^{q-1}\gamma_i U_N^{n+i}}^2  +  \frac{C(\dt)^2 N^2}{\kappa}  \sum_{i=0}^{q-1} \norm{U_N^{n+i}}^2 \\
		{}& - \frac{\beta\dt}{\varepsilon}  \brac{\int A_{M}(W_N^{n+2}, \dots, W_N^{n+q}) - \int A_{M}(W_N^{n+1}, \dots, W_N^{n+q-1}) + \norm{\sum_{i=1}^q c_i W_N^{n+i}}^2_M} .
	\end{aligned}
\end{equation}

Set
\begin{equation}\nonumber
	G_{A_0, U}^n = \int G_{A_0}(U_N^{n}, \dots, U_N^{n+q-1}), \quad A_{M, W}^n = \int A_{M}(W_N^{n+1}, \dots, W_N^{n+q-1})
\end{equation}
and
\begin{equation}\nonumber
	E^n = G_{A_0, U}^{n} + \frac{\beta\dt}{\varepsilon} A_{M, W}^{n}.
\end{equation}
Note that 
\begin{equation} \label{stability-err}
    C^{-1} \sum_{i=0}^{q-1} \norm{U_N^{n+i}}^2 \leq G_{A_0, U}^{n}\leq C \sum_{i=0}^{q-1} \norm{U_N^{n+i}}^2, \quad 0\leq  A_{M, W}^{n} \leq C \sum_{i=1}^{q-1} \norm{U_N^{n+i}}^2.
\end{equation}
By taking $\kappa=d_1/2$ and $\dt\le c_{CFL}/N^2$, it follows from \eqref{eq:combine-estimate} that
\begin{equation}\nonumber
	E^{n+1} - E^n \le C\dt  \sum_{i=0}^{q-1} \norm{U_N^{n+i}}^2 \le C\dt E^n.
\end{equation}
Therefore, we have 
\begin{equation}\nonumber
	E^{n+1} \le (1+C\dt)E^n
\end{equation}
and furthermore 
\begin{equation}\nonumber
	E^n \le e^{CT}E^0.
\end{equation}
Hence we have
\begin{equation}\nonumber
	\norm{U_N^n}^2\le C\sum_{i=0}^{q-1}\brac{\norm{U_N^i}^2 + \frac{\dt}{\varepsilon}\norm{W_N^i}^2}
\end{equation}
and the proof is complete.

\end{proof}

\subsection{Regularity} \label{subsec:regularity-estimate}

To analyze the truncation error of the IMEX-BDF scheme \eqref{scheme1}, we need the following uniform-in-$\varepsilon$ regularity estimate.

For this purpose, we multiply the both sides of \eqref{eq:PDE} with $U^TA_0$ and integrate over $x$ to obtain  
\begin{equation}\nonumber
	\begin{aligned}
		\int U^TA_0U_t + \int U^TA_0  A U_{x} = \frac{1}{\varepsilon} \int U^TA_0QU 
	\end{aligned}
\end{equation}
and thereby
\begin{equation}\label{regu-inequ}
	\begin{aligned}
		\half \int (U^TA_0U)_t + \half \int  (U^TA_0A U)_{x} = \frac{1}{\varepsilon} \int W^TA_{02}\hat{S}W  \le 0.
	\end{aligned}
\end{equation}
Due to the periodic boundary conditions, we have
\begin{equation}\nonumber
	\int U^T(x,t)A_0U(x,t) dx \le \int U^T(x,0)A_0U(x,0) dx,
\end{equation}
which implies
\begin{equation}\nonumber
	\norm{U(\cdot,t)}\le C \norm{U(\cdot,0)}, \quad t \geq 0.
\end{equation}
Here $C$ only depends on the symmetrizer $A_0$. Since \eqref{eq:PDE} is linear with constant coefficients, the partial derivative $\partial_x^s U$ of order $s$ also satisfies \eqref{eq:PDE} and therefore 
\begin{equation}\label{equ:hs}
	\norm{U(\cdot,t)}_{H^s} \le C \norm{U(\cdot,0)}_{H^s}.
\end{equation}
Here $\norm{U(\cdot,t)}_{H^s}$ denotes the standard  norm for the Sobolev space $H^s$ of the periodic function $U = U(x,t)$.

\begin{thm}\label{thm:regularity-const}
	For any integer $s\ge0$, the solution to \eqref{eq:PDE} satisfies
	\begin{enumerate}
		\item 
		for all $t\ge 0$,
		\begin{equation}\label{eq:regularity-Hs}
			\norm{U(\cdot,t)}_{H^s}^2 \le C \norm{U(\cdot,0)}_{H^s}^2,
		\end{equation}
		
		\item
		for all $t\ge  2\delta_0^{-1} s\varepsilon\log(1/\varepsilon)$,
		\begin{equation}\label{eq:regularity-U}
			\norm{\partial_t^{r_1}\partial_x^{r_2}U(\cdot,t)}^2 \le C \norm{U(\cdot,0)}_{H^s}^2, \quad r_1+{r_2}\le s
		\end{equation}
		and
		\begin{equation}\label{eq:regularity-W}
			\norm{\partial_t^{r_1}\partial_x^{r_2}W(\cdot,t)}^2 \le C \varepsilon^2 \norm{U(\cdot,0)}_{H^s}^2, \quad r_1+{r_2}\le s-1.
		\end{equation}	
  Here $\delta_0>0$ is a constant determined by the SPD matricies $A_{02}$ and $A_{02}\hat{S}$, $C$ is a generic constant independent of $\varepsilon$, $U=\begin{pmatrix}
      V\\W
  \end{pmatrix}$, and $r_1, r_2$ are non-negative integers.
	\end{enumerate}
\end{thm}
\begin{proof}
	Estimate \eqref{eq:regularity-Hs} is just \eqref{equ:hs} and
\eqref{eq:regularity-W} simply follows from \eqref{eq:regularity-U} together with the equation
	\begin{equation}\nonumber
		W = \varepsilon \hat{S}^{-1}(W_t + A_{21}V_x + A_{22}W_x).
	\end{equation}
 
 Next, we prove \eqref{eq:regularity-U} by induction on $s$. It is trivial for $s=0$. Assume \eqref{eq:regularity-U} for $(s-1)$ and we prove the estimate with $s$.
 Notice that for any $0\le r\le s-1$, $\partial_t\partial_x^r U$ satisfies the same equation \eqref{eq:PDE}. 
 As in obtaining \eqref{regu-inequ}, we have
	\begin{equation}\nonumber
		\begin{aligned}
			\half  \int ((\partial_t\partial_x^rU)^TA_0 \partial_t\partial_x^rU)_t ={}& \frac{1}{\varepsilon} \int (\partial_t\partial_x^rW)^TA_{02}\hat{S} \partial_t\partial_x^rW \\
    \leq{}& - \frac{\delta_0}{2\varepsilon} \int (\partial_t\partial_x^rW)^T A_{02} \partial_t\partial_x^rW\\
			\le{}& -\frac{\delta_0}{2\varepsilon}\int (\partial_t\partial_x^rU)^T A_0 \partial_t\partial_x^rU + \frac{\delta_0}{2\varepsilon}\int (\partial_t\partial_x^rV)^T A_{01} \partial_t\partial_x^rV.
		\end{aligned}
	\end{equation}
 Here $\delta_0>0$ is a constant determined by the SPD matricies $A_{02}$ and $A_{02}\hat{S}$.

 	Denote
	\begin{equation}\nonumber
		E(t) = \int (\partial_t\partial_x^rU)^TA_0 \partial_t\partial_x^rU.
	\end{equation}
	The last inequality can be written as 
	\begin{equation}\nonumber
		E'(t) \le -\frac{\delta_0}{\varepsilon}E(t) + \frac{C \delta_0}{\varepsilon}\norm{\partial_t\partial_x^rV(t)}^2.
	\end{equation}
By Gronwall's inequality, we have 
	\begin{equation}\label{equ:E-estimate}
		\begin{aligned}
			E(t) \le{}& e^{-\frac{\delta_0}{\varepsilon}t} E(0) + \frac{C\delta_0}{\varepsilon}\int_0^t e^{\frac{\delta_0}{\varepsilon}(\tau-t)}\norm{\partial_t\partial_x^rV(\tau)}^2 d\tau.
		\end{aligned}
	\end{equation}

On the other hand, from the equation for $W$ in \eqref{eq:PDE}
	\begin{equation}\nonumber
		\begin{aligned}
			\partial_t\partial_x^r W = -\partial_x^r(A_{21}V_x + A_{22}W_x - \frac{1}{\varepsilon}\hat{S}W) = - A_{21}\partial_x^{r+1}V - A_{22}\partial_x^{r+1}W + \frac{1}{\varepsilon}\hat{S}\partial_x^r W
		\end{aligned}
	\end{equation}
and estimate \eqref{eq:regularity-Hs}, it follows that 
	\begin{equation}\nonumber
		\norm{\partial_t\partial_x^r W}^2\le C (\frac{1}{\varepsilon^2}+1)\norm{U(\cdot,0)}_{H^s}^2.
	\end{equation}	
	Similarly, we have 
	\begin{equation}\nonumber
		\norm{\partial_t\partial_x^r V}^2\le C \norm{U(\cdot,0)}_{H^s}^2.
	\end{equation}
Thus, it follows from \eqref{equ:E-estimate} that
	\begin{equation}\nonumber
		\begin{aligned}
			E(t) 
			\le{}& Ce^{-\frac{\delta_0}{\varepsilon}t}(\frac{1}{\varepsilon^2}+1)\norm{U(\cdot,0)}_{H^s}^2 + C(1-e^{-\frac{\delta_0}{\varepsilon}t})\norm{U(\cdot,0)}_{H^s}^2 \\
			\le{}& C (\frac{1}{\varepsilon^2}e^{-\frac{\delta_0}{\varepsilon}t}+1)\norm{U(\cdot,0)}_{H^s}^2.
		\end{aligned}
	\end{equation}
 Here we have used
 \begin{equation*}
     E(0) \leq C \norm{\partial_t\partial_x^r V(0)}^2 + C \norm{\partial_t\partial_x^r W(0)}^2 \leq C(\frac{1}{\varepsilon^2}+1) \norm{U(\cdot,0)}_{H^s}^2.
 \end{equation*}
	Then for $t_0=2 \delta_0^{-1}\varepsilon\log(1/\varepsilon)$, we have $E(t_0)\le C\norm{U(\cdot,0)}_{H^s}^2$ and thus
	\begin{equation}\nonumber
		\norm{\partial_t\partial_x^r U(t_0)}^2\le C\norm{U(\cdot,0)}_{H^s}^2.
	\end{equation}

	Now define $\tilde{U}(t)=\partial_tU(t+t_0)$, then $\tilde{U}$ also satisfies the same equation and
	\begin{equation}\nonumber
		\norm{\tilde{U}(0)}_{H^{s-1}}^2 \le C\norm{U(\cdot,0)}_{H^s}^2.
	\end{equation}
	By the induction hypothesis
	\begin{equation}\nonumber
		\norm{\partial_t^{r_1}\partial_x^{r_2} \tilde{U}(t)}^2\le C\norm{U(\cdot,0)}_{H^s}^2, \quad r_1 + r_2 \le s-1, \quad t \ge 2\delta_0^{-1}(s-1)\varepsilon\log(1/\varepsilon),
	\end{equation}
	which implies \eqref{eq:regularity-U}.

\end{proof}


\subsection{Error estimates}\label{subsec:uniform-accuracy}

In this subsection, we establish our main result on the uniform-in-$\varepsilon$ accuracy of the IMEX-BDF scheme \eqref{scheme1}.
As in \cite{hu_uniform_2019}, we consider two types of initial conditions.
\begin{itemize}
	\item 
	\textbf{Type 1}: The initial data $U(x, 0)$ satisfies
	\begin{equation}\nonumber		\norm{\partial_t^{q+1}U(\cdot, 0)}_{H^1} + \norm{\partial_t^{q}U(\cdot, 0)}_{H^2} \le C
	\end{equation}
 for $q=1,2,3,4$. Such data will be used for the IMEX-BDF scheme starting at $T_0\ge0$.

	\item 
	\textbf{Type 2}: The initial data $U(x, 0)$ satisfies
	\begin{equation}\nonumber
		\norm{U(\cdot, 0)}_{H^{q+2}} \le C
	\end{equation}
 for $q=1,2,3,4$. Such data will be used for the IMEX-BDF scheme starting at
 $T_0\ge 2\delta_0^{-1} (q+2)\varepsilon\log(1/\varepsilon)$.
\end{itemize}

\begin{lem}\label{lemma-truncation}
   Let $U^{n}=U(x, t^n)$ is an exact solution to equation \eqref{eq:PDE} with period initial data $U = U(x, 0)$ above. Then the truncation error of the IMEX-BDF \eqref{scheme1} satisfies 
   \begin{equation}\nonumber
       \norm{\sum_{i=0}^q\alpha_i U^{n+i} + \dt A \sum_{i=0}^{q-1}\gamma_i (U^{n+i})_x - \frac{\beta\dt}{\varepsilon} QU^{n+q}}\leq C(\dt)^{q+1}.
   \end{equation}   
\end{lem}

\begin{proof}
Notice that $\partial_t^{r_1}\partial_x^{r_2}U$ satisfies the equation \eqref{eq:PDE}. The regularity estimate \eqref{eq:regularity-Hs} implies 
	\begin{equation}\nonumber
	\norm{\partial_t^{q+1}U(t)}_{H^1} \le C \norm{\partial_t^{q+1}U(0)}_{H^1} \le C
	\end{equation}
	and
	\begin{equation}\nonumber
		\norm{\partial_t^{q}\partial_x U(t)}_{H^1} \le \norm{\partial_t^{q} U(t)}_{H^2} \le C \norm{\partial_t^{q} U(0)}_{H^2}\le C,
	\end{equation}
	for initial data of Type 1.
For initial data of Type 2, the regularity estimate \eqref{eq:regularity-U} leads to
	\begin{equation}\nonumber		\norm{\partial_t^{r_1}\partial_x^{r_2}U(t)}\le C \norm{U(t)}_{H^{q+2}} \le C, \quad r_1+r_2\le q+2,
	\end{equation}
	for any $t\ge 2\delta_0^{-1}(q+2)\varepsilon\log(1/\varepsilon)$.
 This implies
	\begin{equation}\nonumber
		\norm{\partial_t^{q+1}U(t)}_{H^1} \le C
	\end{equation}
 by taking $r_1=q+1$, $r_2=0, 1$ and 
	\begin{equation}\nonumber
		\norm{\partial_t^{q}\partial_xU(t)}_{H^1} \le C
	\end{equation}
 by taking $r_1=q$, $r_2=0, 1, 2$.
	Moreover, it follows from the Sobolev inequality that 
	\begin{equation}\nonumber
		\norm{\partial_t^{q+1}U(t)}_{L^{\infty}} +\norm{\partial_t^{q}\partial_xU(t)}_{L^{\infty}} \le C\left(
		\norm{\partial_t^{q+1}U(t)}_{H^1} + \norm{\partial_t^{q}\partial_xU(t)}_{H^1} \right )\le C
	\end{equation}
 for the initial data of the two types.

On the other hand, from \cite{hundsdorfe2007imex} we know the following facts related to the IMEX-BDF scheme:
	\begin{equation}\nonumber
		\abs{\sum_{i=0}^q \alpha_i u^{n+i} - \beta\dt\partial_t u^{n+q}} \le C\dt^{q+1}\max_{t\in[T^0,T]}\abs{\partial_t^{q+1}u^{n+q}}
	\end{equation}
	and
	\begin{equation}\nonumber
		\abs{\sum_{i=0}^{q-1}\gamma_i\partial_x u^{n+i} - \beta\dt\partial_x u^{n+q}} \le C\dt^{q+1}\max_{t\in[T^0,T]}\abs{\partial_t^{q}\partial_xu^{n+q}}
	\end{equation}
 for any smooth function $u=u(x, t)$, where $u^n:=u(x,t^n)$. Thus, for the spatially  periodic function $u=u(x, t)$, we have
	\begin{equation}\nonumber
		\norm{\sum_{i=0}^q \alpha_i u^{n+i} - \beta\dt\partial_t u^{n+q}} \le C\norm{\sum_{i=0}^q \alpha_i u^{n+i} - \beta\dt\partial_t u^{n+q}}_{L^{\infty}} \le C(\dt)^{q+1}
	\end{equation}
	and
	\begin{equation}\nonumber
		\norm{\sum_{i=0}^{q-1}\gamma_i\partial_x u^{n+i} - \beta\dt\partial_x u^{n+q}} \le C\norm{\sum_{i=0}^{q-1}\gamma_i\partial_x u^{n+i} - \beta\dt\partial_x u^{n+q}}_{L^{\infty}} \le C(\dt)^{q+1}.
	\end{equation}

	Denote by $R_U^n$ the truncation error of the IMEX-BDF scheme \eqref{scheme1}:
	\begin{equation}\nonumber
		R_U^n = \sum_{i=0}^q\alpha_i U^{n+i} + \dt A \sum_{i=0}^{q-1}\gamma_i (U^{n+i})_x - \frac{\beta\dt}{\varepsilon} QU^{n+q}.
	\end{equation}
 It follows from the last two inequalities that
	\begin{equation}\nonumber
		\begin{aligned}
		    \norm{R_U^n}& = \norm{\sum_{i=0}^q\alpha_i U^{n+i} + \dt A \sum_{i=0}^{q-1}\gamma_i (U^{n+i})_x - \frac{\beta\dt}{\varepsilon} QU^{n+q}}\\[4mm]
      &\leq  \norm{\sum_{i=0}^q \alpha_i U^{n+i} - \beta\dt\partial_t U^{n+q}} + \norm{\beta\dt\partial_t U^{n+q} + \beta\dt A \partial_x U^{n+q} - \frac{\beta\dt}{\varepsilon} QU^{n+q} }\\[4mm]
      &~~  + \norm{\dt A \sum_{i=0}^{q-1}\gamma_i (U^{n+i})_x - \beta\dt A \partial_x U^{n+q}}\\
      &\leq C (\dt)^{q+1} + 0 +  C (\dt)^{q+1}\\
      & \le C(\dt)^{q+1}.
		\end{aligned}
	\end{equation}
 This completes the proof.
\end{proof}

\begin{thm}\label{thm:uniform-accuracy}
Under the conditions of Theorem \ref{thm:stability-const}, the IMEX-BDF scheme \eqref{scheme1} for system \eqref{eq:PDE} is uniformly $q$-th order accurate, that is
\begin{equation}\nonumber
	\norm{U(\cdot, t_n) - U_N^n}^2 \le C\left( (\dt)^{2q} + e_{init} \right).
\end{equation}
Here $U=U(x, t)$ is the exact solution to equation \eqref{eq:PDE} with initial data above,
 $C$ is a constant independent of $\varepsilon$, $N$ and $\dt$, and $e_{init}$ is related to the initial projection error
\begin{equation}\nonumber
	e_{init} := \sum_{i=0}^{q-1}\left( \norm{V(\cdot, t_i) - V_N^i}^2 + (1+\frac{\dt}{\varepsilon})\norm{W(\cdot, t_i) - W_N^i} \right)
\end{equation}
with 
 $U=\begin{pmatrix}
      V\\W
  \end{pmatrix}$ and $U^n_N=\begin{pmatrix}
      V^n_N\\W^n_N
  \end{pmatrix}$.
\end{thm}

\begin{proof}
	Set $\delta U^n=U(x, t_n) - U_N^n$. It is clear that the error $\delta U^n$ satisfies the scheme \eqref{scheme1}  with residue $R_U^n$. Then by repeating the argument of Theorem \ref{thm:stability-const} and using
 	\begin{equation}\nonumber
		\begin{aligned}
			\int \left( \delta U^{n+q} - \sum_{i=1}^q\eta_i\delta U^{n+i} \right) R_U^n dx 
			\le{}& \kappa\dt \norm{\delta U^{n+q}}^2 + \kappa C\dt\sum_{i=1}^{q-1} \norm{\delta U^{n+i}}^2 + \frac{C}{\kappa}(\dt)^{2q+1}
		\end{aligned}
	\end{equation}
	with $\kappa>0$, we obtain
	\begin{equation}\label{err-sccuracy}
		\begin{aligned}
		    E^{n+1} - E^n \le& C \dt  \sum_{i=0}^{q-1} \norm{\delta U^{n+i}}^2 - \frac{\beta\dt}{\varepsilon}   \norm{\sum_{i=1}^q c_i \delta W^{n+i}}^2_M + \kappa\dt\norm{\delta U^{n+q}}^2 + \frac{C}{\kappa}(\dt)^{2q+1}\\
		\end{aligned}
	\end{equation}
	where
 	\begin{equation}\nonumber
		E^n = G_{A_0, \delta U}^{n} + \frac{\beta\dt}{\varepsilon} A_{M, \delta W}^{n}
	\end{equation}
 with
	\begin{equation}\nonumber
		G_{A_0, \delta U}^n = \int G_{A_0}(\delta U^{n}, \dots, \delta U^{n+q-1}), \quad A_{M, \delta W}^n = \int A_{M}(\delta W^{n+1}, \dots, \delta W^{n+q-1}).
	\end{equation}
As inequalities \eqref{stability-err}, we have
\begin{equation} \nonumber
    C^{-1} \sum_{i=0}^{q-1} \norm{\delta U^{n+i}}^2 \leq G_{A_0, \delta U}^{n}\leq C \sum_{i=0}^{q-1} \norm{\delta U^{n+i}}^2, \quad 0 \leq A_{M, \delta W}^{n}  \leq C \sum_{i=1}^{q-1} \norm{\delta U^{n+i}}^2.
\end{equation}
Then inequality \eqref{err-sccuracy} gives 
	\begin{equation}\nonumber
		\begin{aligned}
		    E^{n+1} - E^n &\leq  C \Delta t G_{A_0, \delta U}^n + C \kappa\dt G_{A_0, \delta U}^{n+1} + \frac{C}{\kappa}(\dt)^{2q+1}\\
      &\leq C \Delta t E^n + C \kappa\dt E^{n+1} + \frac{C}{\kappa}(\dt)^{2q+1}.
		\end{aligned}
	\end{equation}
	With this, we take $\kappa$ sufficiently small to obtain 
	\begin{equation}\nonumber
		E^{n+1} \le (1+C\dt)E^n + C(\dt)^{2q+1}
	\end{equation}
	implying 
	\begin{equation}\nonumber
		E^n \le C E^0 + C(\dt)^{2q}.
	\end{equation}
Since
\begin{equation*}
    \norm{U(\cdot, t_n) - U_N^n}^2 = \norm{\delta U^n}^2 \leq C E^n
\end{equation*}
and
\begin{equation*}
    E^0 \leq  C \sum_{i=0}^{q-1} \left( \norm{\delta V^{i}}^2 + (1+\frac{\dt}{\varepsilon})\norm{\delta W^{i}}^2 \right) =  C e_{init},
\end{equation*}
the proof is completed.
\end{proof}

We end this section with
the following corollary.
\begin{cor}
	Under the conditions of Theorem \ref{thm:uniform-accuracy}, 
 if 
	\begin{equation}\nonumber
		\norm{U(\cdot, T_0)}_{H^{2q+1}} + \norm{\partial_t U(\cdot, T_0)}_{H^{2q}} \le C,
	\end{equation}
	the error estimate 
	\begin{equation}\nonumber
		\norm{U(\cdot, t_n) - U_N^n}^2 \le C \left( (\dt)^{2q} + \frac{1}{N^{4q}} \right),
	\end{equation}
	holds for integer $n$ such that $t_n=T_0+n\dt \le T$.
\end{cor}
\begin{proof}
By Theorem \ref{thm:uniform-accuracy}, it suffices to prove that $e_{init}\le C/N^{4q}$. To do so,
    we use the following
property of Fourier projection \cite{hesthaven2007spectral} and Theorem \ref{thm:regularity-const} to obtain
	\begin{equation}\nonumber
		\norm{U(\cdot, t_i) - U_N^i}^2 \le \frac{1}{N^{4q+2}}\norm{U(\cdot, t_i)}^2_{H^{2q+1}} \le \frac{C}{N^{4q+2}}.
	\end{equation}
	Similarly, we have
	\begin{equation}\nonumber
		\norm{\partial_x (V(\cdot, t_i) - V_N^i)}^2 + \norm{\partial_t (W(\cdot, t_i) - W_N^i)}^2 \le \frac{C}{N^{4q}}.
	\end{equation}
 Then we deduce from the equation for $W$ in \eqref{eq:PDE} that
	\begin{equation}\nonumber
		\norm{W(\cdot, t_i) - W_N^i}^2 \le \varepsilon^2 \frac{C}{N^{4q}}.
	\end{equation}
 Hence $e_{init}\le C/N^{4q}$ and the conclusion follows.
\end{proof}

\section{Numerical tests}\label{sec:numerical-tests}

In this section, we numerically test the accuracy of the IMEX-BDF schemes applied to several linearized hyperbolic relaxation systems including the Aw-Rascle-Zhang traffic model \cite{aw2000siam,zhang2002non}, the Broadwell model \cite{broadwell1964shock}, and the Grad's moment system \cite{grad1949kinetic,cai2014cpam}.
In all the numerical tests, we adopt the Fourier-Galerkin spectral method for spatial discretization with modes $|k|\leq N$ and fix $N=100$ to ensure that the discretization error in space is much smaller than that in time.
The reference solution $U_{ref}$ is computed with a much finer time step.

\subsection{Aw-Rascle-Zhang traffic model}
The model 
\cite{aw2000siam,zhang2002non} is 
\begin{equation}\nonumber
    \begin{aligned}
        &\partial_t \rho + \partial_x (\rho v)  = 0,\\
        &\partial_t v  + \left(v - \rho p'(\rho) \right)\partial_x v = \frac{V(\rho) - v}{\varepsilon},
    \end{aligned}
\end{equation}
with
\begin{equation}\nonumber
   p(\rho) = c_0 \rho^\gamma, \quad V(\rho) = v_f \left(1 - \frac{\rho}{\rho_m}\right). 
\end{equation}
Here $\rho = \rho(x, t)$ is the traffic density, $v = v(x,t)$ is the traffic speed, and $\varepsilon$ is a relaxation time characterizing the response of the drivers to the traffic situation. 
The variable $p(\rho)$ is the traffic pressure and the equilibrium velocity-density relationship $V(\rho)$ is given in the Greenshield model \cite{greenshields1935study}.
The linearization of the model  around a uniform steady state $(\rho^\star, v^\star)$ is
\begin{equation}\nonumber
    \begin{aligned}
        &\partial_t \rho + v^\star \partial_x \rho + \rho^\star \partial_x v  = 0,\\
        &\partial_t v  - \left(\rho^\star p'(\rho^\star) - v^\star \right)\partial_x v = \frac{\rho V'(\rho^\star) - v}{\varepsilon}.
    \end{aligned}
\end{equation}

In our numerical test, we take
\begin{equation}\nonumber
    c_0 =\frac{3}{2}, \quad \gamma = 1, \quad \rho_m = 8, \quad v_f = 4, \quad (\rho^\star, v^\star) = (1, 1).
\end{equation}
Then the linearized model becomes
\begin{equation}\nonumber
    \begin{aligned}
        \partial_t U + A\partial_x U = \frac{1}{\varepsilon} Q U,
    \end{aligned} 
\end{equation}
with 
\begin{equation}\nonumber
    \begin{aligned}
       U = (\rho, ~ v)^T, \quad A = \begin{pmatrix}
           1 & 1\\
           0 &-\frac{1}{2}
       \end{pmatrix},\quad Q=\begin{pmatrix}
           0 & 0\\
           -\frac{1}{2} & -1
       \end{pmatrix}.
    \end{aligned} 
\end{equation}
It is easy to verify that the last system satisfies the structural stability condition with 
\begin{equation*}
    \begin{aligned}
        P = \begin{pmatrix}
            1 & 0\\
            \frac{1}{2} & 1
        \end{pmatrix}, \quad A_0 = \begin{pmatrix}
           3 & 2\\
           2 & 4
        \end{pmatrix}. 
    \end{aligned}
\end{equation*}

The computational domain is $[0, 1]$ with periodic boundary conditions and the initial data are given by 
\begin{equation}\nonumber
    \rho(x,0) = \sin(2\pi x) + 1.1.
\end{equation}
For the second-order scheme, we choose 
the initial data for $v$ as
\begin{equation}\nonumber
    v(x, 0) = -\frac{1}{2}\rho(x, 0),
\end{equation}
which is consistent up to $O(1)$. For the third-order scheme, we choose 
the initial data for $v$ as
\begin{equation}\nonumber
    v(x, 0) = -\frac{1}{2}\rho(x, 0) -  \frac{\varepsilon}{2} \partial_x \rho(x, 0),
\end{equation}
which is consistent up to $O(\varepsilon)$.
For the fourth-order scheme, we choose 
the initial data for $v$ as
\begin{equation}\nonumber
    v(x, 0) = -\frac{1}{2}\rho(x, 0) -  \frac{\varepsilon}{2} \partial_x \rho(x, 0) - \frac{\varepsilon^2}{4} \partial_{xx} \rho(x, 0),
\end{equation}
which is consistent up to $O(\varepsilon^2)$.
The starting values of the IMEX-BDF scheme at $t = i\Delta t$ with $i =1, \cdots, q-1$, are prepared using the IMEX-RK schemes (ARS(2,2,2) for second- and third-order scheme, ARS(4,4,3) for fourth-order scheme \cite{ascher1997}) 
with a much smaller time step $\delta t = \Delta t/500$. We compute the solution to time $T=1$ and estimate the $L^2$ error of the solutions $U_{\Delta t}$ as $\norm{U_{\Delta t}-U_{ref}}$.


Table \ref{tab:ARZ-error} gives the $L^2$ error and convergence rates with respect to $\Delta t$ of IMEX-BDF schemes of order $q=2, 3, 4$ with $\varepsilon$ ranging from $10^{-7}$ to $1$. We can observe that the numerical results are in perfect agreement with our theoretical analysis for various values of $\varepsilon$. The minor order degeneration in the fourth-order scheme with $\Delta t=1.79\times10^{-4}$ is due to the machine precision limitations.

\begin{table}[htbp]
  \centering
  \caption{Aw-Rascle-Zhang traffic model: The $L^2$ error of the solutions computed by IMEX-BDF schemes of order $q=2,3,4$.}
  \label{tab:ARZ-error}%
    \begin{tabular}{c|c|c|c|c|c|c|c}
    \hline
    \multirow{2}[4]{*}{$\varepsilon$} & \multirow{2}[4]{*}{$\Delta t$} & \multicolumn{2}{c|}{second order} & \multicolumn{2}{c|}{third order} & \multicolumn{2}{c}{fourth order} \bigstrut\\
    \cline{3-8}  &  &   $L^2$-error & order  & $L^2$-error & order  & $L^2$-error & order \bigstrut\\   
    \hline
    \multirow{4}[3]{3em}{$10^{-7}$}
        & 1.43e-03 & 4.46e-04 & -    & 2.25e-06 & -    & 1.08e-08 & -    \\
        & 7.14e-04 & 1.11e-04 & 2.00 & 2.82e-07 & 3.00 & 6.74e-10 & 4.00 \\
        & 3.57e-04 & 2.75e-05 & 2.02 & 3.52e-08 & 3.00 & 4.24e-11 & 3.99 \\
        & 1.79e-04 & 6.55e-06 & 2.07 & 4.34e-09 & 3.02 & 3.14e-12 & 3.76 \\
    \hline
    \multirow{4}[3]{3em}{$10^{-6}$}
        & 1.43e-03 & 4.46e-04 & -    & 2.25e-06 & -    & 1.08e-08 & -    \\
        & 7.14e-04 & 1.11e-04 & 2.00 & 2.82e-07 & 3.00 & 6.74e-10 & 4.00 \\
        & 3.57e-04 & 2.75e-05 & 2.02 & 3.52e-08 & 3.00 & 4.23e-11 & 4.00 \\
        & 1.79e-04 & 6.55e-06 & 2.07 & 4.34e-09 & 3.02 & 3.13e-12 & 3.76 \\
    \hline
    \multirow{4}[3]{3em}{$10^{-5}$}
        & 1.43e-03 & 4.46e-04 & -    & 2.25e-06 & -    & 1.08e-08 & -    \\
        & 7.14e-04 & 1.11e-04 & 2.00 & 2.82e-07 & 3.00 & 6.74e-10 & 4.00 \\
        & 3.57e-04 & 2.75e-05 & 2.02 & 3.52e-08 & 3.00 & 4.24e-11 & 3.99 \\
        & 1.79e-04 & 6.56e-06 & 2.07 & 4.34e-09 & 3.02 & 3.14e-12 & 3.76 \\
    \hline
    \multirow{4}[3]{3em}{$10^{-4}$}
        & 1.43e-03 & 4.46e-04 & -    & 2.25e-06 & -    & 1.07e-08 & -    \\
        & 7.14e-04 & 1.12e-04 & 1.99 & 2.81e-07 & 3.00 & 6.73e-10 & 4.00 \\
        & 3.57e-04 & 2.83e-05 & 1.98 & 3.52e-08 & 3.00 & 4.22e-11 & 3.99 \\
        & 1.79e-04 & 7.11e-06 & 1.99 & 4.40e-09 & 3.00 & 2.97e-12 & 3.83 \\
    \hline
    \multirow{4}[3]{3em}{$10^{-3}$}
        & 1.43e-03 & 4.98e-04 & -    & 2.29e-06 & -    & 1.05e-08 & -    \\
        & 7.14e-04 & 1.34e-04 & 1.89 & 3.05e-07 & 2.91 & 6.76e-10 & 3.96 \\
        & 3.57e-04 & 3.52e-05 & 1.93 & 4.09e-08 & 2.90 & 4.46e-11 & 3.92 \\
        & 1.79e-04 & 8.67e-06 & 2.02 & 5.31e-09 & 2.95 & 3.01e-12 & 3.89 \\
   \hline
    \multirow{4}[3]{3em}{$10^{-2}$}
        & 1.43e-03 & 5.81e-04 & -    & 2.71e-06 & -    & 7.40e-09 & -    \\
        & 7.14e-04 & 1.47e-04 & 1.98 & 3.48e-07 & 2.96 & 4.69e-10 & 3.98 \\
        & 3.57e-04 & 3.67e-05 & 2.00 & 4.41e-08 & 2.98 & 2.89e-11 & 4.02 \\
        & 1.79e-04 & 8.77e-06 & 2.06 & 5.48e-09 & 3.01 & 1.26e-12 & 4.52 \\
    \hline
    \multirow{4}[3]{3em}{$10^{-1}$}
        & 1.43e-03 & 4.07e-04 & -    & 2.26e-06 & -    & 1.72e-08 & -    \\
        & 7.14e-04 & 9.97e-05 & 2.03 & 2.83e-07 & 3.00 & 1.08e-09 & 4.00 \\
        & 3.57e-04 & 2.46e-05 & 2.02 & 3.53e-08 & 3.00 & 6.71e-11 & 4.00 \\
        & 1.79e-04 & 5.86e-06 & 2.07 & 4.36e-09 & 3.02 & 4.15e-12 & 4.01 \\
    \hline
    \multirow{4}[3]{3em}{$10^{0}$}
        & 1.43e-03 & 2.77e-03 & -    & 4.32e-05 & -    & 1.27e-06 & -    \\
        & 7.14e-04 & 3.37e-04 & 3.04 & 5.41e-06 & 3.00 & 7.94e-08 & 4.00 \\
        & 3.57e-04 & 8.34e-05 & 2.02 & 6.76e-07 & 3.00 & 4.97e-09 & 4.00 \\
        & 1.79e-04 & 1.97e-05 & 2.07 & 8.34e-08 & 3.02 & 3.19e-10 & 3.96 \\                     
    \hline   
    \end{tabular}%
\end{table}%


\subsection{Broadwell model}

The Broadwell model is a simplified discrete velocity model for the Boltzmann equation \cite{broadwell1964shock}. It describes a two-dimensional (2D) gas as composed of particles of only four velocities with a binary collision law and spatial variation in only one direction. When looking for one-dimensional solutions of the 2D gas, the evolution equations of the model are given by 
\begin{equation}\nonumber
	\begin{aligned}
		\partial_t f_{+} + \partial_x f_{+} &= -\frac{1}{\varepsilon}(f_{+}f_{-} - f_0^2),\\
		\partial_t f_{-} - \partial_x f_{-} &= -\frac{1}{\varepsilon}(f_{+}f_{-} - f_0^2),\\
		\partial_t f_{0} &= \frac{1}{\varepsilon}(f_{+}f_{-} - f_0^2).
	\end{aligned}
\end{equation}
Here $f_{+}$, $f_{-}$ and $f_0$ denote the particle density function at time $t$, position $x$ with velocity $1$, $-1$ and $0$, respectively, $\varepsilon>0$ is the mean free path. 
Set 
\begin{equation}\nonumber
	\rho = f_{+} + 2f_0 + f_{-}, \quad m = f_{+} - f_{-}, \quad z = f_{+} + f_{-}.
\end{equation}
The Broadwell equations can be rewritten as 
\begin{equation}\nonumber
	\begin{aligned}
		\partial_t \rho + \partial_x m &= 0,\\
		\partial_t m + \partial_x z &= 0,\\
		\partial_t z + \partial_x m &=  \frac{1}{2\varepsilon}( \rho^2 + m^2 - 2\rho z).
	\end{aligned}
\end{equation}

A local Maxwellian is the density function that satisfies $z = \frac{1}{2\rho}(\rho^2 + m^2)$. Considering the linearized version at $\rho_\star = 2,  m_\star=0, z_\star=1 $, we obtain the linearized Broadwell system as follows
\begin{equation}\nonumber
    \begin{aligned}
        \partial_t U + A \partial_x U = \frac{1}{\varepsilon}QU,
    \end{aligned}
\end{equation}
with 
\begin{equation}\nonumber
    \begin{aligned}
        U =  (\rho, ~m, ~z)^T, \qquad A=\begin{pmatrix}
            0 & 1 & 0\\
            0 & 0 & 1\\
            0 & 1 & 0\\
        \end{pmatrix}, \qquad Q=\begin{pmatrix}
            0 & 0 & 0\\
            0 & 0 & 0\\
            1 & 0 & -2\\
            \end{pmatrix}.
    \end{aligned}
\end{equation}
It has been shown in \cite{yong_singular_1999} that the Broadwell model satisfies the structural stability condition.

In our numerical test, the computational spatial domain is $[-\pi, \pi]$ with periodic boundary conditions and the initial data of $\rho$ and $m$ are given by
\begin{equation}\nonumber
    \begin{aligned}
        \rho(x, 0) = 1 + a_{\rho}\sin(2x), \quad  m(x, 0) = \rho(x, 0) \left( \frac{1}{2} + a_{u}\cos(2x) \right)\\
    \end{aligned}
\end{equation}
with $a_\rho = 0.3$ and $a_u =0.1$.
For the second-order scheme, we choose the initial data for $z$ as
\begin{equation}\nonumber
    z(x, 0) = \frac{1}{2}\rho(x, 0),
\end{equation}
which is consistent up to $O(1)$. For the third-order scheme, we choose the initial data for $z$ as
\begin{equation}\nonumber
    z(x, 0) = \frac{1}{2}\rho(x,0) - \frac{\varepsilon}{4}\partial_x m(x, 0),
\end{equation}
which is consistent up to $O(\varepsilon)$.
 For the fourth-order scheme, we choose the initial data for $z$ as
\begin{equation}\nonumber
    z(x, 0) = \frac{1}{2}\rho(x,0) - \frac{\varepsilon}{4}\partial_x m(x, 0)  - \frac{\varepsilon^2}{16}\partial_{xx} \rho(x, 0),
\end{equation}
which is consistent up to $O(\varepsilon^2)$.
The starting values at $i\Delta t$ with $i =1, \cdots, q-1$ are prepared using ARS(4,4,3) with a much smaller time step $\delta t = \Delta/500$. We compute the solution to time $T=2$ and estimate the error of the solutions $U_{\Delta t}$ as $\norm{U_{\Delta t}-U_{ref}}$.



In Table \ref{tab:Broadwel-error}, we present the numerical results of IMEX-BDF schemes of order $q = 2, 3, 4$, and various values of $\Delta t$ and $\varepsilon$. The uniform $q-$th order accuracy is clearly achieved for $q = 2, 3, 4$. This closely aligns with our theoretical analysis.



\begin{table}[htbp]
  \centering
  \caption{Broadwell system: The $L^2$ error of the solutions computed by IMEX-BDF schemes of order $q=2,3,4$.}
  \label{tab:Broadwel-error}%
    \begin{tabular}{c|c|c|c|c|c|c|c}
    \hline
    \multirow{2}[4]{*}{$\varepsilon$} & \multirow{2}[4]{*}{$\Delta t$} & \multicolumn{2}{c|}{second order} & \multicolumn{2}{c|}{third order} & \multicolumn{2}{c}{fourth order} \bigstrut\\
    \cline{3-8}  &  &   $L^2$-error & order  & $L^2$-error & order  & $L^2$-error & order \bigstrut\\   
    \hline
    \multirow{4}[3]{3em}{$10^{-7}$}
        & 5.00e-03 & 4.58e-04 & -    & 4.07e-06 & -    & 5.03e-08 & -    \\
        & 2.50e-03 & 1.14e-05 & 2.00 & 5.09e-07 & 3.00 & 3.17e-09 & 3.99 \\
        & 1.25e-03 & 2.83e-05 & 2.02 & 6.36e-08 & 3.00 & 1.98e-10 & 4.00 \\
        & 6.25e-04 & 6.74e-06 & 2.07 & 7.84e-09 & 3.02 & 1.24e-11 & 4.00 \\
    \hline
    \multirow{4}[3]{3em}{$10^{-6}$}
        & 5.00e-03 & 4.58e-04 & -    & 4.07e-06 & -    & 5.03e-08 & -    \\
        & 2.50e-03 & 1.14e-05 & 2.00 & 5.09e-07 & 3.00 & 3.17e-09 & 3.99 \\
        & 1.25e-03 & 2.83e-05 & 2.02 & 6.36e-08 & 3.00 & 1.99e-10 & 4.00 \\
        & 6.25e-04 & 6.74e-06 & 2.07 & 7.84e-09 & 3.02 & 1.25e-11 & 4.00 \\
    \hline
    \multirow{4}[3]{3em}{$10^{-5}$}
        & 5.00e-03 & 4.59e-04 & -    & 4.07e-06 & -    & 5.03e-08 & -    \\
        & 2.50e-03 & 1.14e-05 & 2.00 & 5.09e-07 & 3.00 & 3.17e-09 & 3.99 \\
        & 1.25e-03 & 2.83e-05 & 2.02 & 6.36e-08 & 3.00 & 1.98e-10 & 4.00 \\
        & 6.25e-04 & 6.75e-06 & 2.07 & 7.84e-09 & 3.02 & 1.24e-11 & 4.00 \\
    \hline
    \multirow{4}[3]{3em}{$10^{-4}$}
        & 5.00e-03 & 4.59e-04 & -    & 4.07e-06 & -    & 5.03e-08 & -    \\
        & 2.50e-03 & 1.15e-05 & 2.00 & 5.09e-07 & 3.00 & 3.17e-09 & 3.99 \\
        & 1.25e-03 & 2.84e-05 & 2.01 & 6.36e-08 & 3.00 & 1.98e-10 & 4.00 \\
        & 6.25e-04 & 6.80e-06 & 2.06 & 7.84e-09 & 3.02 & 1.25e-11 & 3.99 \\
    \hline
    \multirow{4}[3]{3em}{$10^{-3}$}
        & 5.00e-03 & 4.64e-04 & -    & 4.06e-06 & -    & 5.02e-08 & -    \\
        & 2.50e-03 & 1.17e-05 & 1.98 & 5.08e-07 & 3.00 & 3.16e-09 & 3.99 \\
        & 1.25e-03 & 2.97e-05 & 1.98 & 6.34e-08 & 2.99 & 1.98e-10 & 4.00 \\
        & 6.25e-04 & 7.28e-06 & 2.06 & 7.83e-09 & 3.01 & 1.23e-11 & 4.00 \\
    \hline
    \multirow{4}[3]{3em}{$10^{-2}$}
        & 5.00e-03 & 4.98e-04 & -    & 3.98e-06 & -    & 4.91e-08 & -    \\
        & 2.50e-03 & 1.27e-05 & 1.97 & 5.01e-07 & 2.99 & 3.09e-09 & 3.99 \\
        & 1.25e-03 & 3.20e-05 & 1.99 & 6.28e-08 & 2.99 & 1.94e-10 & 4.00 \\
        & 6.25e-04 & 7.68e-06 & 2.06 & 7.77e-09 & 3.01 & 1.20e-11 & 4.00 \\
    \hline
    \multirow{4}[3]{3em}{$10^{-1}$}
        & 5.00e-03 & 5.00e-04 & -    & 3.69e-06 & -    & 4.04e-08 & -    \\
        & 2.50e-03 & 1.25e-05 & 2.00 & 4.64e-07 & 2.99 & 2.53e-09 & 3.99 \\
        & 1.25e-03 & 3.11e-05 & 2.01 & 5.81e-08 & 3.00 & 1.59e-10 & 4.00 \\
        & 6.25e-04 & 7.41e-06 & 2.07 & 7.17e-09 & 3.02 & 9.80e-12 & 4.02 \\
    \hline  
    \multirow{4}[3]{3em}{$10^{0}$}
        & 5.00e-03 & 5.13e-04 & -    & 7.61e-06 & -    & 1.17e-07 & -    \\
        & 2.50e-03 & 1.28e-05 & 2.00 & 9.51e-07 & 3.00 & 7.31e-09 & 4.00 \\
        & 1.25e-03 & 3.16e-05 & 2.02 & 1.19e-07 & 3.00 & 4.57e-10 & 4.00 \\
        & 6.25e-04 & 7.52e-06 & 2.07 & 1.46e-08 & 3.02 & 2.85e-11 & 4.00 \\
    \hline   
    \end{tabular}%
\end{table}%


\subsection{Linearized Grad's moment system} 
The linearized Grad's moment system in 1D \cite{grad1949kinetic,cai2014cpam,zhao2017stability} reads as
\begin{equation}\label{equ:linear-moment-equ}
	\partial_t U + A\partial_x U = \frac{1}{\varepsilon}QU
\end{equation}
with
\begin{equation}\nonumber
	\begin{aligned}
		U = \begin{pmatrix}
			\rho \\ w \\ \theta/\sqrt{2} \\ \sqrt{3!}f_3 \\ \vdots \\\sqrt{M!}f_M
		\end{pmatrix},
		A = \begin{pmatrix}
			0 & 1 & & & &  \\
			1 & 0 & \sqrt{2} & & & \\
			& \sqrt{2} & 0 & \sqrt{3} & &   \\
			& & \sqrt{3} & 0 & \ddots &   \\
			& & & \ddots & 0 & \sqrt{M} \\
			& & & & \sqrt{M} & 0			
		\end{pmatrix},
		Q = -\diag(0, 0, 0, \underbrace{1, \cdots, 1}_{M-2}).
	\end{aligned}
\end{equation}
In the above equation, $\rho$ is the density, $w$ is the macroscopic velocity, $\theta$ is the temperature and $f_3, \cdots, f_M$ with $M\geq 3$ are high order moments. The moment system is obtained by taking moments on the both sides of the Bhatnagar-Gross-Krook (BGK) model \cite{Bhatnagar1954511}. 
It was shown in \cite{Di2017nm,zhao2017stability} that the moment system satisfies the structural stability condition.
Here we only consider its linearized version. 

The spatial domain is taken as $x\in [-\pi, \pi]$ with periodic boundary conditions.
We solve the linearized Grad’s moment system \eqref{equ:linear-moment-equ} with $M=5$.  
The initial data are prepared by
\begin{equation}\nonumber
	(\rho, ~ w, ~ \theta)(x, 0) = \left( \sin(2x)+1.1, ~ 0, ~ \sqrt{2} \right),\qquad (f_3, f_4, f_5) = (0, ~ 0, ~ 0).
\end{equation}
The starting values at $i\Delta t, i = 0, \cdots, q-1$, are prepared using an ARS(4,4,3) scheme with a much smaller time step $\delta t = \Delta t/500$. We compute the solution to time $T=1$ and estimate the error of the solution $U_{\Delta t}$ as $\norm{U_{\Delta t}-U_{ref}}$.

Table \ref{tab:LinearMoment-error} show the $L^2$ error of IMEX-BDF schemes of order $q=2, 3, 4$, and various values of $\Delta t$ and $\varepsilon$. Again, we observe the uniform accuracy of the scheme with $\varepsilon$ ranging from $10^{-7}$ to 1.


\begin{table}[htbp]
  \centering
  \caption{Linearized Grad's moment system: The $L^2$ error of the solutions computed by IMEX-BDF schemes.}
  \label{tab:LinearMoment-error}%
    \begin{tabular}{c|c|c|c|c|c|c|c}
    \hline
    \multirow{2}[4]{*}{$\varepsilon$} & \multirow{2}[4]{*}{$\Delta t$} & \multicolumn{2}{c|}{second order} & \multicolumn{2}{c|}{third order} & \multicolumn{2}{c}{fourth order} \bigstrut\\
    \cline{3-8}  &  &   $L^2$-error & order  & $L^2$-error & order  & $L^2$-error & order \bigstrut\\   
    \hline
    \multirow{4}[3]{3em}{$10^{-7}$}
        & 2.50e-03 & 1.04e-03 & -    & 1.34e-05 & -    & 9.29e-08 & -    \\
        & 1.25e-03 & 2.62e-04 & 2.00 & 1.68e-06 & 2.99 & 5.88e-09 & 3.98 \\
        & 6.25e-04 & 6.49e-05 & 2.01 & 2.10e-07 & 3.00 & 3.70e-10 & 3.99 \\
        & 3.13e-04 & 1.55e-05 & 2.07 & 2.59e-08 & 3.02 & 2.34e-11 & 3.98 \\
    \hline
    \multirow{4}[3]{3em}{$10^{-6}$}
        & 2.50e-03 & 1.04e-03 & -    & 1.34e-05 & -    & 9.29e-08 & -    \\
        & 1.25e-03 & 2.62e-04 & 2.00 & 1.68e-06 & 2.99 & 5.88e-09 & 3.98 \\
        & 6.25e-04 & 6.49e-05 & 2.01 & 2.10e-07 & 3.00 & 3.70e-10 & 3.99 \\
        & 3.13e-04 & 1.55e-05 & 2.07 & 2.59e-08 & 3.02 & 2.33e-11 & 3.99 \\
    \hline
    \multirow{4}[3]{3em}{$10^{-5}$}
        & 2.50e-03 & 1.04e-03 & -    & 1.34e-05 & -    & 9.29e-08 & -    \\
        & 1.25e-03 & 2.62e-04 & 2.00 & 1.68e-06 & 2.99 & 5.88e-09 & 3.98 \\
        & 6.25e-04 & 6.49e-05 & 2.01 & 2.10e-07 & 3.00 & 3.70e-10 & 3.99 \\
        & 3.13e-04 & 1.55e-05 & 2.07 & 2.59e-08 & 3.02 & 2.31e-11 & 4.00 \\
    \hline
    \multirow{4}[3]{3em}{$10^{-4}$}
        & 2.50e-03 & 1.04e-03 & -    & 1.34e-05 & -    & 9.29e-08 & -    \\
        & 1.25e-03 & 2.62e-04 & 2.00 & 1.68e-06 & 2.99 & 5.88e-09 & 3.98 \\
        & 6.25e-04 & 6.49e-05 & 2.01 & 2.10e-07 & 3.00 & 3.70e-10 & 3.99 \\
        & 3.13e-04 & 1.55e-05 & 2.07 & 2.59e-08 & 3.02 & 2.31e-11 & 4.00 \\
    \hline
    \multirow{4}[3]{3em}{$10^{-3}$}
        & 2.50e-03 & 1.04e-03 & -    & 1.34e-05 & -    & 9.31e-08 & -    \\
        & 1.25e-03 & 2.62e-04 & 2.00 & 1.68e-06 & 2.99 & 5.92e-09 & 3.97 \\
        & 6.25e-04 & 6.49e-05 & 2.01 & 2.10e-07 & 3.00 & 3.75e-10 & 3.98 \\
        & 3.13e-04 & 1.55e-05 & 2.07 & 2.59e-08 & 3.02 & 2.36e-11 & 3.99 \\
    \hline
    \multirow{4}[3]{3em}{$10^{-2}$}
        & 2.50e-03 & 1.05e-03 & -    & 1.33e-05 & -    & 9.58e-08 & -    \\
        & 1.25e-03 & 2.62e-04 & 2.00 & 1.67e-06 & 2.99 & 6.09e-09 & 3.98 \\
        & 6.25e-04 & 6.50e-05 & 2.01 & 2.09e-07 & 3.00 & 3.84e-10 & 3.99 \\
        & 3.13e-04 & 1.55e-05 & 2.07 & 2.58e-08 & 3.02 & 2.41e-11 & 4.00 \\
    \hline
    \multirow{4}[3]{3em}{$10^{-1}$}
        & 2.50e-03 & 1.06e-03 & -    & 1.31e-05 & -    & 1.01e-07 & -    \\
        & 1.25e-03 & 2.66e-04 & 2.00 & 1.65e-06 & 2.99 & 6.40e-09 & 3.98 \\
        & 6.25e-04 & 6.59e-05 & 2.01 & 2.06e-07 & 3.00 & 4.02e-10 & 3.99 \\
        & 3.13e-04 & 1.57e-05 & 2.07 & 2.54e-08 & 3.02 & 2.52e-11 & 4.00 \\
    \hline  
    \multirow{4}[3]{3em}{$10^{0}$}
        & 2.50e-03 & 1.07e-03 & -    & 1.26e-05 & -    & 1.01e-07 & -    \\
        & 1.25e-03 & 2.68e-04 & 2.00 & 1.58e-06 & 3.00 & 6.39e-09 & 3.99 \\
        & 6.25e-04 & 6.64e-05 & 2.01 & 1.98e-07 & 3.00 & 4.01e-10 & 3.99 \\
        & 3.13e-04 & 1.58e-05 & 2.07 & 2.44e-08 & 3.02 & 2.51e-11 & 4.00 \\
    \hline   
    \end{tabular}%
\end{table}%

\bibliographystyle{abbrv}

\end{document}